\documentclass[11pt]{amsart}

\usepackage{epsfig}
\usepackage{amssymb}
\usepackage{enumerate}

\usepackage{a4wide}

\def\cqfd{\mbox{}\nolinebreak\hfill$\Box$\medbreak\par}

\newcommand{\Z}{\mathbb{Z}}

\newcommand{\B}{\mathbb{B}}
\renewcommand{\P}{\mathbb{P}}

\newcommand{\Tau}{\mathcal{T}}

\theoremstyle{plain}

\newtheorem{thm}{Theorem}[section]

\newtheorem{lem}[thm]{Lemma}
\newtheorem{prop}[thm]{Proposition}

\theoremstyle{definition}
\newtheorem{dfn}[thm]{Definition}
\newtheorem{rk}[thm]{Remark}

\numberwithin{equation}{section}

\newenvironment{pf}{\noindent\textbf{Proof:}}{\cqfd}

\title{Forgetful maps between Deligne-Mostow ball quotients}

\author{Martin Deraux}

\begin{document}

\date{\today}

\address{\textsc
  Universit\'e de Grenoble I\\
  Institut Fourier, UMR 5582\\
  BP 74\\
  38402 Saint-Martin d'H\`eres\\
  France\\
}
\email{deraux@fourier.ujf-grenoble.fr} 

\begin{abstract}
  We study forgetful maps between Deligne-Mostow moduli spaces of
  weighted points on $\P^1$, and classify the forgetful maps that
  extend to a map of orbifolds between the stable completions. The
  cases where this happens include the Livn{\'e} fibrations and the
  Mostow/Toledo maps between complex hyperbolic surfaces. They also
  include a retraction of a $3$-dimensional ball quotient onto one of
  its $1$-dimensional totally geodesic complex submanifolds.
\end{abstract}

\maketitle

\section{Introduction}

The goal of this paper is to collect some information about known maps
between Deligne-Mostow ball quotients of various dimensions. Only one
map given here is new (it gives a non-trivial map from a
$3$-dimensional ball quotient to a compact Riemann surface), but we
also find it worthwhile to commit to print the fact that forgetful
maps (in the sense of the present paper, see
section~\ref{sec:forgetful}) cannot yield any other examples.

The simplest examples of holomorphic maps between ball quotients are
given by unbranched coverings, obtained simply by taking subgroups of
finite index of the fundamental group of the relevant ball quotient.
Another class of examples is given by totally geodesic maps, which are
also easily constructed between ball quotients of any dimension. In
fact, it is well known that there are many holomorphic totally
geodesic inclusions between Deligne-Mostow quotients, and we only
briefly review how to describe those maps in
Proposition~\ref{prop:inclusions}.

Another way to obtain non-trivial holomorphic maps is to construct
branched coverings. This can of course easily be done for Riemann
surfaces, which appear in the present paper in the form of surjective
homomorphisms between various triangle groups. In higher dimensions
however, simple branched coverings of ball quotients cannot be ball
quotients themselves (see the computation of characteristic classes
that comes up in the Mostow-Siu construction,
see~\cite{mostowsiu},~\cite{derauxmathann} and~\cite{deraux3fold}),
and it is not clear how to determine whether a given ball quotient
admits branched coverings (simple or not). Note that the maps between
complex hyperbolic surfaces constructed by Mostow and Toledo
(see~\cite{toledomaps}) are certainly not simple branched coverings,
in fact they branch around certain complex totally geodesic curves,
but they also contract some such curves.

Finally, one might hope to get certain ball quotients to fiber over
ball quotients, and this was first achieved by Livn{\'e} in his
thesis, see~\cite{livne}. His ball quotient are actually closely
related to some of the Deligne-Mostow lattices, see \S16
of~\cite{delignemostowbook}, and the corresponding fibrations can then
be interpreted in terms of the forgetful map construction presented in
this paper; this remark was the basis for the construction of maps to
Riemann surfaces used in the author's thesis
(see~\cite{derauxmathann}). As mentioned above, in this paper, we show
that the forgetful map construction can be generalized to give a
similar fibration of a $3$-dimensional compact ball quotient to a
Riemann surface, but we do not know of any method to obtain fibrations
for higher-dimensional examples.

The main result of this paper can be thought as shedding some light on
the general question of existence of surjective holomorphic map
$X^m\rightarrow Y^n$ between compact ball quotients when $m>n\geq 2$
(see the question raised by Siu in~\cite{siuconjecture}, p.~182). One
might think of the statement of~Theorem~\ref{thm:main} as giving some
evidence for the non-existence of such maps, but note that the ball
quotients considered here are not particularly representative (there
are only finitely many Deligne-Mostow ball quotients in dimension
$\geq 2$), and we are only considering a very specific construction of
maps between them.

Our results should also be put in perspective with a recent result of
Koziarz and Mok, see~\cite{koziarzmok}, that precludes the existence
of \emph{submersive} holomorphic maps between ball quotients. The maps
obtained in this paper are indeed not submersive, and some explicit
fibers are in fact singular divisors.

It should be pointed out that the ball with its Bergman metric is a
\emph{rank one} Hermitian symmetric space, which makes it very
different from irreducible higher rank Hermitian symmetric spaces.
Indeed, if $X=\Gamma\setminus\Omega$ is a compact manifold modelled on
a higher rank bounded symmetric domain $\Omega$, then there is no
surjective holomorphic map from $X$ onto any nonpositively curved
Hermitian manifold apart from geodesic coverings
(see~\cite{mokannals}).

Another interesting feature of some of the maps that appear in this
paper is that they exhibit a retraction of the relevant ball quotient
onto one of its totally geodesic submanifolds. In the context of
\emph{real} hyperbolic geometry, such retractions have been obtained for
certain arithmetic real hyperbolic manifolds (see~\cite{bergeron}),
without any restriction on the dimension.

This paper was written as an answer to various questions asked over
the years by Domingo Toledo, Sai-Kee Yeung and Ngaiming Mok, whom the
author wishes to thank for their interest in this work. The existence
of a map from a $3$-ball quotient to a compact hyperbolic Riemann
surface as in Theorem~\ref{thm:main}(\ref{it:3dim}) was also known to
Sai-Kee Yeung.

\section{Review of Deligne-Mostow theory}

\subsection{The Picard integrality condition}\label{sec:int}

We start by collecting some facts from Deligne-Mostow theory
(see~\cite{delignemostow},~\cite{mostowihes}), following the
exposition in~\cite{klw}. We state only what is needed for the
purpose of this paper (for a more thorough survey see~\cite{looijenga}
for instance).

Given an integer $m\geq 1$, we would like to consider various
structures on the moduli space of $m+3$ points on $\P^1$. In order to
form a geometric invariant theory quotient of $(\P^1)^{m+3}$, we need
to pick a line bundle $\mathcal{L}$ on $(\P^1)^{m+3}$ and a lift of
the $PGL_2$-action to $\mathcal{L}$ (we refer to this data as a
polarization). We shall choose various polarizations, each encoded by
a choice of weights $\mu=(\mu_1,\dots,\mu_{m+3})\in]0,1[^{m+3}$ for
the $(m+3)$ points; throughout the paper, the weights shall be taken
to be rational numbers, and we assume moreover that
$\sum_{j=1}^{m+3}\mu_j=2$.

The line bundle $\mathcal{L}_\mu$ on $(\P^1)^{m+3}$ associated to
$\mu$ is given by $\displaystyle\underset{j}{\boxtimes}\
\mathcal{O}(2d\mu_j)$ where $d$ is the common denominator of the
$\mu_j$ (see section~4.6 of~\cite{delignemostow}). The corresponding
geometric invariant theory quotient has a simple description, as we
now recall.

We define $M$ to be the set of $(m+3)$-tuples of pairwise distinct
points on $\P^1$, and the following chain of subsets
$$
M\subset M^{\mu}_{st} \subset M^{\mu}_{sst}
$$
of $(\P^1)^{m+3}$ by allowing only certain coincidences of points.
$M^{\mu}_{st}$ denotes the subset of $(m+3)$-tuples where we allow
coincidence of points only when the sum of the corresponding weights
is $<1$. The set $M^{\mu}_{sst}$ is defined similarly, allowing
coincidence of points whose weights add up to $\leq 1$.

For each strictly semistable point $x$, there is a unique partition
$\{1,\dots,m+3\}=S_1\sqcup S_2$ such that for some $j=1$ or $2$, the
points $x_i$ with indices $i\in S_j$ coincide. We then define the
corresponding quotient spaces
$$
Q\subset Q^\mu_{st} \subset Q^\mu_{sst}
$$
where $Q_{st}^\mu$ is the set of $PGL_2$-orbits of points of $M_{st}$,
and two strictly semistable points $x$ and $y$ are identified if and
only if the associated partitions coincide.

Note that the space $Q_{sst}^\mu$ is compact, but in general it is
singular (whereas $Q_{st}^\mu$ is always smooth).

For convenience, we will sometimes write $Q^\mu$ for $Q$, and $M^\mu$
for $M$, even though these spaces depend on $\mu$ only through the
number $m+3$ of components of $\mu$.

\begin{dfn}
  We denote by $D_{ij}^\mu$ the image in $Q^\mu_{st}$ of the set of
  points $(x_1,\dots,x_{m+3})\in M^\mu_{st}$ with $x_i=x_j$. 
\end{dfn}
When the dependence on $\mu$ is clear, we sometimes write $D_{ij}$ for
$D_{ij}^\mu$. This set is a divisor in $Q_{st}^\mu$ only if
$\mu_i+\mu_j<1$.

\begin{dfn}
The set of weights $\mu$ is said to satisfy the Picard integrality
condition if
$$
(1-\mu_i-\mu_j)^{-1}\in\Z\qquad \textrm{(INT)}
$$
whenever $i\neq j$ and $\mu_i+\mu_j<1$.

For any $i\neq j\in\{1,\dots,m+3\}$, we shall write
\begin{equation}
  \label{eq:orbifoldWeights}
  d_{ij}^{(\mu)}=(1-\mu_i-\mu_j)^{-1}
\end{equation}
When the Picard integrality condition holds, $d^{(\mu)}_{ij}$ is
always an integer (or infinity), regardless of whether or not
$\mu_i+\mu_j<1$, see~\cite{delignemostow}, page~26. When no confusion
arises, we shall simply write $d_{ij}$ for $d_{ij}^{(\mu)}$.

\end{dfn}
When the Picard integrality condition is satisfied, the main result
of~\cite{delignemostow} gives $Q_{st}^\mu$ the structure of a complex
hyperbolic orbifold in the following sense.
\begin{thm}\label{thm:dm}
  If $\mu$ satisfies condition INT, then there is a lattice
  $\Gamma_\mu$ in $PU(m,1)$ such that the orbifold 
  $$
  X_\mu=\Gamma_\mu\setminus \B^m
  $$
  has an underlying smooth complex manifold structure isomorphic to
  $Q_{st}^\mu$. Under this identification, the singular locus of
  $X_\mu$ consists of the points of $Q_{st}^\mu - Q$, and the divisor
  $D_{ij}^\mu$ has weight $d_{ij}^{(\mu)}$.
\end{thm}

In other words, for every torsion-free subgroup $G_\mu\subset
\Gamma_\mu$ of finite index, the map $G_\mu\setminus \B^m\rightarrow
\Gamma_\mu\setminus\B^m$ can be thought of as giving a description of
$G_\mu\setminus \B^m$ as a branched covering of $Q_{st}^\mu$, with
ramification of index $d_{ij}^{(\mu)}$ above $D_{ij}^\mu$. 

The orbifold fundamental group of $X_\mu=\Gamma_\mu\setminus \B^m$ is
of course just $\Gamma_\mu$, and an explicit presentation for that
group can be deduced from Lemma~\ref{lem:kernel} below.  We follow the
notation used in~\cite{delignemostow} and denote by
$\widetilde{Q^{\mu}}$ the preimage in the ball of $Q^{\mu}$. Since the
map $\widetilde{Q^{\mu}}\rightarrow Q^{\mu}$ is an unbranched
covering, the fundamental group $\pi_1(\widetilde{Q^{\mu}})$ is
identified with a subgroup $K_\mu$ of $\pi_1(Q^{\mu})$.

For reasons that are explained in~\cite{delignemostow}, we shall call
$\widetilde{Q^{\mu}}$ the \emph{monodromy cover}, since its
fundamental group is in fact the kernel of the monodromy
representation (see~\cite{delignemostow}, section~8).

The following result is a special case of Lemma~8.6.1
in~\cite{delignemostowbook}.
\begin{lem}\label{lem:kernel}
  Let $\gamma_{ij}$ be a small loop that goes once around
  $D_{ij}^\mu$, i.e. a loop that corresponds to $x_i$ turning once
  around $x_j$, see Figure~\ref{fig:twist}~(left). Then
  $$\Gamma_\mu\simeq \pi_1(Q^\mu)/K_\mu$$
  where $K_\mu$ is the normal subgroup of $\pi_1(Q^\mu)$ generated by
  the $\gamma_{ij}^{d_{ij}}$, $\mu_i+\mu_j<1$. 
\end{lem}

\begin{rk}
  The lattice $\Gamma_\mu$ is cocompact if and only if
  $Q_{st}^\mu=Q_{sst}^\mu$, i.e. no subset of the weights adds up to
  exactly $1$. When $Q_{st}^\mu$ is not compact and the corresponding
  lattice $\Gamma_\mu$ is arithmetic, the compactification
  $Q_{sst}^\mu$ is homeomorphic to the Baily-Borel compactification of
  $\Gamma_\mu\setminus \B^m$ (in the non-arithmetic cases, the
  compactification is obtained by adding a finite number of
  cusps).
\end{rk}

A nice feature of the above picture is that the divisors $D_{ij}^\mu$
themselves have a modular interpretation, as moduli spaces of $m+2$
points on $\P^1$, with two of the weights $\mu_i$ and $\mu_j$ replaced
by their sum.

More generally, the configurations obtained by letting certain subsets
of the $m+3$ points coalesce give suborbifolds of larger codimension,
and they also have a modular interpretation (we only allow points to
coalesce if the sum of the corresponding weights is strictly less than
$1$).

\begin{dfn}
  For any subset $I\subset \{1,\dots,m+3\}$ consisting of $r+1$
  elements, the \emph{contraction} of $\mu=(\mu_1,\dots,\mu_{m+3})$
  along $I$ is the $(n-r)$-tuple obtained by replacing the weights
  $\mu_j$, $j\in I$ by their sum, i.e. 
  $\mu^{(I)}=(\mu_{i_1},\dots,\mu_{i_{n-r-1}},\sum_{i\in I}\mu_i$),
  where $\{1,\dots,n\}\setminus I = \{i_1,\dots,i_{n-r-1}\}$.
\end{dfn}

We shall consider the contraction $\mu^{(I)}$ only when $\sum_{i\in
  I}\mu_i<1$, in which case $\mu^{(I)}$ satisfies the running
hypotheses of this section, hence defines a lattice
$\Gamma_{\mu^{(I)}}$ acting on a ball of dimension $m-r$ as in
Theorem~\ref{thm:dm}.

\begin{dfn}
  $\mu^{(I)}$ is called a \emph{hyperbolic contraction} of $\mu$ if
  $\sum_{i\in I}\mu_i<1$.
\end{dfn}

\begin{prop} \label{prop:inclusions} Let $\mu$ satisfy condition INT,
  and let $\mu^{(I)}$ be a hyperbolic contraction of $\mu$.
  Then $\mu^{(I)}$ also satisfies INT, and moreover there exists a
  totally geodesic subball $B\in \B^{m}$ of codimension $r=|I|-1$ such
  that the image of $B$ in $\Gamma_\mu\setminus \B^{m}$ is isomorphic
  as orbifolds to $\Gamma_{\mu^{(I)}}\setminus \B^{m-r}$.
  $\Gamma_{\mu^{(I)}}$ is isomorphic to the stabilizer of $B$ modulo
  its fixed point stabilizer.
\end{prop}

This proposition follows from~(8.8.1) in~\cite{delignemostow},
see also Lemma~2.4 in~\cite{mostowdiscontinuous}. It gives totally
geodesic inclusions between various Deligne-Mostow orbifolds.

\subsection{Condition $\frac{1}{2}$INT} \label{sec:halfint}

We now discuss how to generalize the results of the previous section
as in~\cite{mostowihes}. The generalized version is also the
one that appears in Thurston's account of this theory,
see~\cite{thurstonshapes}. The idea is to consider moduli of $(m+3)$
\emph{unordered} points rather than ordered. Since we consider
weighted points, we allow identification of $(m+3)$-tuples of points
that differ by ordering only when the corresponding permutation of the
indices preserves the weights.

More specifically, we fix a partition $m+3=i_1+\dots+i_k$, and
consider moduli of sets $S_1\dots,S_k$ of points on $\P^1$, with $S_j$
having cardinality $i_j$ for each $j=1,\dots,k$. In terms of the
notation used in the previous section, this moduli space is a quotient
of $Q^\mu/\Sigma$, where
$\Sigma=\Sigma_{i_1}\times\dots\times\Sigma_{i_k}$ is a product of
symmetric groups.

In order to describe the choice of weights in this setting in terms of
the notation used in the previous section, we consider $(m+3)$-tuples
$\mu$ with $0<\mu_j<1$ for all $j$ and $\sum \mu_j=2$, and break up
the index set $I=\{1,\dots,m+3\}$ as a disjoint union
$I_1\sqcup\dots\sqcup I_k$ in such a way that, for each $j$, the
$\mu_i$, $i\in I_j$ are equal. In the sequel we shall always assume
that the index sets $I_j$ are arranged in increasing order, in the
sense that if $j<j'$, all the elements of $I_j$ are smaller than those
of $I_{j'}$.  Note that, by construction, $\mu$ is then invariant
under $\Sigma$.
\begin{dfn}
  The pair $\mu,\Sigma$ as above satisfies the half-integrality
  condition $\frac{1}{2}$INT if for all $i\neq j$ such that
  $\mu_i+\mu_j<1$, we have
  
  $$
  (1-\mu_i-\mu_j)^{-1}\in\left\{\begin{array}{ll}
      \Z & \textrm{ if $i$ and $j$ are not in the same $\Sigma$-orbit}\\
      \frac{1}{2}\Z & \textrm{ if $i$ and $j$ are in the same
        $\Sigma$-orbit}
      \end{array}\right.
  $$
\end{dfn}
We adapt the definition of the integers $d_{ij}$ accordingly, and
set
\begin{equation}\label{eq:sigmaintweights}
d_{ij}^{(\mu,\Sigma)}=\left\{\begin{array}{ll}
    (1-\mu_i-\mu_j)^{-1}\textrm{ if $i$ and $j$ are not in the same $\Sigma$-orbit}\\
    2(1-\mu_i-\mu_j)^{-1}\textrm{ if $i$ and $j$ are in the same
      $\Sigma$-orbit}
      \end{array}\right.
\end{equation}
and write $D_{ij}^{\mu,\Sigma}$ for the image of $D_{ij}^\mu\subset
Q_{st}^\mu$ in the quotient $Q_{st}^\mu/\Sigma$. As above, when no
confusion arises, we shall simply write $d_{ij}$ instead of
$d_{ij}^{(\mu,\Sigma)}$.

\begin{rk}
  \begin{enumerate}
  \item We do not necessarily assume that $\mu_i\neq \mu_{i'}$ when
    $i\in I_j$ and $i'\in I_{j'}$ with $j\neq j'$. In other words, we
    do not assume that the sets of indices $I_j$ are as large as
    possible to get $\mu$ to be $\Sigma$-invariant.
  \item When $i$ and $j$ are in the same $\Sigma$-orbit, we do not
    assume that $(1-\mu_i-\mu_j)^{-1}$ are in $\frac{1}{2}\Z\setminus\Z$.
  \end{enumerate}
\end{rk}

The action of $\Sigma$ on $M$ clearly descends to an action on $Q$,
and this action extends to an action on $Q_{st}^\mu$. In general, the
quotient space $Q_{st}^\mu/\Sigma$ has singularities, but the content
of the main result of~\cite{mostowihes} is that it carries a complex
hyperbolic orbifold structure, similar to the one mentioned in the
previous section:
\begin{thm}
  \label{thm:mostow}
  If $\mu,\Sigma$ satisfies condition $\frac{1}{2}$INT, then there is
  a lattice $\Gamma_{\mu,\Sigma}$ in $PU(m,1)$ such that the orbifold
  $$
  X_{\mu,\Sigma}=\Gamma_{\mu,\Sigma}\setminus \B^m
  $$ 
  has the same underlying (singular) algebraic variety as
  $Q_{st}^\mu/\Sigma$. Under this identification, the divisors
  $D_{ij}^{\mu,\Sigma}$ have weight $d_{ij}^{(\mu,\Sigma)}$, and the
  other divisors with weight $>1$ have weight two, and are the images
  of codimension one fixed point sets of elements of $\Sigma$ that are
  contained in $Q^\mu$.
\end{thm}

\begin{rk}\label{rk:fixedpoints}
  \begin{enumerate}
  \item There are indeed sometimes elements of $\Sigma$ that fix a
    codimension one subset contained in $Q^\mu$. The list of cases
    where that happens can be deduced from in Lemma~8.3.2
    of~\cite{delignemostowbook} (the elements that give codimension
    one fixed point set contained in $Q^\mu$ are bitranspositions when
    $m+3=5$, and tritranspositions when $m+3=6$).
  \item The same criterion as in the previous section determines
    whether the relevant ball quotient is compact or not, namely
    $\Gamma_{\mu,\Sigma}$ is cocompact if and only if no subset of the
    weights adds up to exactly $1$.
  \end{enumerate}
\end{rk}

\begin{figure}
  \centering
  \epsfig{figure=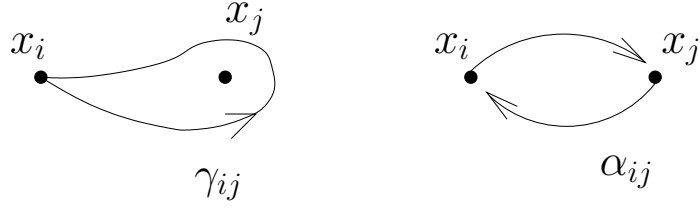, width=0.6\textwidth}
  \caption{A full twist (left) and a half twist (right) between $x_i$
    and $x_j$. $\gamma_{ij}$ induces a loop in $Q^\mu_{st}$ that goes
    once around $D_{ij}^\mu$. When $i$ and $j$ are not in the same
    $\Sigma$-orbit, $\gamma_{ij}$ induces a loop in
    $Q^\mu_{st}/\Sigma$ that goes once around $D_{ij}^{\mu,\Sigma}$.
    When $i$ and $j$ are in the same $\Sigma$-orbit, $\alpha_{ij}$
    induces a loop in $Q_{st}^\mu/\Sigma$ that goes once around
    $D_{ij}^{\mu,\Sigma}$, and $\alpha_{ij}^2$ induces the same loop
    as $\gamma_{ij}$.}
  \label{fig:twist}
\end{figure}

The analogue of Lemma~\ref{lem:kernel} in the context of
$\frac{1}{2}$INT examples is slightly more complicated to state.
\begin{dfn}\label{dfn:freeaction}
  Let $Q'^{\mu,\Sigma}$ denote the largest open set of $Q^\mu$ on
  which the action of $\Sigma$ is free.
\end{dfn}
When $i$ and $j$ are not in the same $\Sigma$-orbit, we use the same
notation as above and write $\gamma_{ij}$ for a full twist between
$x_i$ and $x_j$. If $i,j$ are in the same $\Sigma$-orbit, we denote by
$\alpha_{ij}$ the corresponding half twist (see
Figure~\ref{fig:twist}). Note that the $\gamma_{ij}$ (resp.
$\alpha_{ij}$) with $\mu_i+\mu_j<1$ give ``small loops'' in
$Q'^{\mu,\Sigma}/\Sigma$
around
$D_{ij}^{\mu,\Sigma}$.

Now consider the elements of $\Sigma$ that have a codimension one
fixed point set contained in $Q^\mu$ (see
Remark~\ref{rk:fixedpoints}(1)), and denote by $B_1,\dots,B_k$ the
components of their image in $Q^{\mu,\Sigma}/\Sigma$. Write $\beta_j$
for a small loop in $Q'^{\mu,\Sigma}/\Sigma$ that goes once around
$B_j$.

\begin{lem}
  \label{lem:kernelhalfint}
  We have
  $$
  \Gamma_{\mu,\Sigma}\simeq
  \pi_1(Q'^{\mu,\Sigma}/\Sigma)/K_{\mu,\Sigma}
  $$
  where 
  $$K_{\mu,\Sigma}=\langle\langle \alpha_{ij}^{d_{ij}}, \beta_i^2\rangle\rangle$$ 
  is the normal subgroup of $\pi_1(Q'^{\mu,\Sigma}/\Sigma)$ generated
  by the $\alpha_{ij}^{d_{ij}}$ such that $\mu_i+\mu_j<1$, and by the
  $\beta_i^2$, $i=1,\dots,k$.
\end{lem}

The group $K_{\mu,\Sigma}$ can again be interpreted as the fundamental
group of a certain unbranched covering
$\widetilde{Q'^{\mu,\Sigma}/\Sigma}$ of $Q'^{\mu,\Sigma}/\Sigma$,
which we call the monodromy cover (see~\cite{mostowihes}). Under the
identification given in Theorem~\ref{thm:mostow},
$\widetilde{Q'^{\mu,\Sigma}/\Sigma}$ identifies with an open set in
the ball $\B^m$.

\subsection{Obvious commensurabilities}

If condition $\frac{1}{2}$INT holds for a given $\mu$ but for two
different symmetry groups, then the corresponding two lattices are
commensurable, as we now explain.

If the partition $I_1,\dots,I_k$ is a refinement of a partition
$J_1,\dots,J_l$ (i.e. every $I_i$ is contained in some $J_j$), and we
denote by $\Sigma^{(I)}$ (resp. $\Sigma^{(J)}$) the corresponding
symmetry group preserving the partition $\{I_i\}$ (resp. preserving
$\{J_j\}$), then clearly $\Sigma^{(I)}$ is a subgroup of
$\Sigma^{(J)}$, hence there is a natural map
$$
Q_{st}^\mu/\Sigma^{(I)}\rightarrow Q_{st}^\mu/\Sigma^{(J)}
$$

Provided that $\mu,\Sigma^{(I)}$ and $\mu,\Sigma^{(J)}$ both satisfy
condition $\frac{1}{2}$INT, one gets accordingly that
$$
\Gamma_{\mu,\Sigma^{(I)}}\subset \Gamma_{\mu,\Sigma^{(J)}}
$$
is a subset of index $[\Sigma^{(J)}:\Sigma^{(I)}]$.

For each $\mu$ that satisfies $\frac{1}{2}$INT for some symmetry
group, there is a \emph{finest} partition of the indices for which
condition $\frac{1}{2}$INT holds, and it is obtained by requiring that
$(1-\mu_{i}+\mu_{j})^{-1}$ be in $\frac{1}{2}\Z$ but \emph{not} in
$\Z$ for each $i\neq j$ in a common subset of the partition.

\begin{rk}\label{rk:duality}
  There is an extra relation between $1$-dimensional Deligne-Mostow
  moduli spaces, alluded to in~\cite{delignemostow}, namely the moduli
  space for $\mu=(\mu_1,\mu_2,\mu_3,\mu_4)$ is isomorphic to the one
  for $(1-\mu_1,1-\mu_2,1-\mu_3,1-\mu_4)$. We shall refer to these two
  sets of weights as being \emph{dual} to each other,
  see~\cite{delignemostow}, p.~84.
\end{rk}

\section{Forgetful maps} \label{sec:forgetful}

The rough idea of our construction is to consider the obvious maps
from the moduli space of $m'=m+3$ points to the moduli space of
$n'=n+3$ points on $\P^1$, obtained by forgetting $m-n$ points (here
we assume that $m\geq n$). This is only a rough idea, because we are
actually interested in moduli spaces of \emph{weighted} points on
$\P^1$.

Without loss of generality, we may assume that the points that get
forgotten are the last $m-n$, i.e. we consider the map
$(\P^1)^{m'}\rightarrow (\P^1)^{n'}$ given by
\begin{equation}
(x_1,\dots,x_{n'},x_{n'+1},\dots,x_{m'})\mapsto (x_1,\dots,x_{n'}).\label{eq:forgetful} 
\end{equation}
This map is clearly $PGL_2$-equivariant, hence it always induces a
holomorphic map $Q^\mu\rightarrow Q^\nu$ for any set of weights
$\mu=(\mu_1,\dots,\mu_{m+3})$ and $\nu=(\nu_1,\dots,\nu_{n+3})$.

We now fix two choices $\mu,\Sigma$ and $\nu,\Tau$ of weights with
symmetry that satisfy condition $\frac{1}{2}$INT, and we wish to
investigate the following question. \emph{When does the above
  forgetful map induce a map of orbifolds between the corresponding
moduli spaces $X_{\mu,\Sigma}$ and $X_{\nu,\Tau}$}?

The fact that the forgetful maps do indeed induce maps of orbifolds
for certain choices of $\mu,\Sigma$ and $\nu,\Tau$ was already noticed
in a couple of places in the literature. For $m=2$ and $n=1$ this was
used in~\cite{derauxthesis} (see also~\cite{derauxmathann}). For
$m=n=2$, it was used by Toledo in~\cite{toledomaps}. Note that in the
equidimensional case, the number of weighted points is the same for
both moduli spaces in the source and target, so the forgetful map is
simply the identity on the level of $Q^{\mu}=Q^{\nu}$ (but it turns
out that its extension to stable moduli spaces contracts some
divisors).

The goal of the present paper is to give the list of all other cases
where the forgetful maps give maps of orbifolds between the orbifold
ball quotients.  For simplicity we focus on the case of compact
orbifolds, i.e. we consider sets of weights $\mu$ such that
$Q^\mu_{st}=Q^\mu_{sst}$ (equivalently that no sum $\sum_{j\in J}
\mu_j$ is equal to one, for any $J\subset\{1,\dots,m+3\}$).

One way to summarize the results of the classification is the
following.

\begin{thm}\label{thm:main}
  Suppose $m\geq n$, $\mu,\Sigma$ and $\nu,\Tau$ satisfy condition
  $\frac{1}{2}$INT, and assume that the orbifolds $X_{\mu,\Sigma}$ and
  $X_{\nu,\Tau}$ are compact, of dimension $m$ and $n$ respectively.
  \begin{enumerate}[(i)]
  \item If $(m,n)$ is not $(1,1)$, $(2,2)$ $(2,1)$ or $(3,1)$, then
    the forgetful map never induces a map of orbifolds.
  \item When $(m,n)=(1,1)$, there are many forgetful maps that induce
    maps of orbifolds (giving among others surjective homomorphisms
    between triangle groups).
  \item When $(m,n)=(2,2)$, the forgetful maps that induce maps of
    orbifolds correspond to the Mostow/Toledo maps.
  \item Many forgetful maps induce maps of orbifolds in the case
    $(m,n)=(2,1)$, some corresponding to the Livn{\'e}
    fibrations.
  \item When $(m,n)=(3,1)$, up to symmetry and obvious
    commensurability (see section~\ref{sec:halfint}), there is
    precisely one forgetful map that yields a map of orbifolds,
    corresponding to the weights $\mu=(3,3,3,3,3,1)/8$ and
    $\nu=(3,3,3,7)/8$.
    \label{it:3dim}
  \end{enumerate}
\end{thm}

As noted in the introduction, Koziarz and Mok have recently shown that
there are no submersive maps $X^m\rightarrow X^n$ between compact ball
quotients apart from unbranched coverings (this holds for finite
volume ball quotients as well provided $n\geq 2$,
see~\cite{koziarzmok}). The map $X^3\rightarrow X^1$ that appears in
part~(\ref{it:3dim}) of the theorem is of course not a submersion (for
a description of the non-submersive locus, see
Remark~\ref{rk:singular}). The fact that the above construction should
produce maps from a compact $3$-ball quotient to a compact hyperbolic
Riemann surface seems not to have appeared anywhere in the literature.

We shall not attempt to use general results from geometric invariant
theory, since the above question can be answered in a fairly
elementary way by using Theorem~\ref{thm:mostow} and
Lemma~\ref{lem:kernelhalfint}. In order for the forgetful map to
induce a map between the orbifold quotients, we shall require that the
forgetful map be compatible with the symmetry groups, in the sense
that
\begin{center}
$Q'^{\mu,\Sigma}$ maps into $Q'^{\nu,\Tau}$, and this map descends to
a map $Q'^{\mu,\Sigma}/\Sigma\rightarrow Q'^{\nu,\Tau}/\Tau$.
\end{center}

In order for the forgetful map to induce a map of orbifolds, we need
to require moreover that the map $Q'^{\mu,\Sigma}/\Sigma\rightarrow
Q'^{\nu,\Tau}/\Tau$ lifts to monodromy covers, which in view of
Lemma~\ref{lem:kernelhalfint} can be equivalently expressed by the
fact that
\begin{center}
  $K_{\mu,\Sigma}$ maps into $K_{\nu,\Tau}$.
\end{center}

If the latter condition holds, then the lift defines a holomorphic map
from the complement of a discrete union of subballs in $\B^m$, which
extends to the whole $\B^{m}$ by Hartogs' theorem.

\subsection{Compatibility of the symmetry groups}

As in the previous paragraphs, we fix two sets of weights with
symmetry $\mu,\Sigma$ and $\nu,\Tau$ that both satisfy condition
$\frac{1}{2}$INT. In order for the map~(\ref{eq:forgetful}) to induce
a map
$$
Q^\mu/\Sigma\rightarrow Q^\nu/\Tau,
$$
we need to require a compatibility condition between the partition
$I=I_1\sqcup\dots\sqcup I_k$ (resp. $J=J_1\sqcup\dots\sqcup J_l$)
corresponding to $\mu,\Sigma$ (resp. to $\nu,\Tau$). 

\begin{lem} \label{lem:compat} The forgetful map $Q^\mu\rightarrow
  Q^\nu$ descends to a map $Q^\mu/\Sigma\rightarrow Q^\nu/\Tau$ if and
  only if for each $i$, $I_i$ is either entirely forgotten (i.e.
  $j>n'$ for all $j\in I_i$), or contained in $J_j$ for some $j$.
\end{lem}

Finally, the condition that $Q'^{\mu,\Sigma}$ be mapped to
$Q'^{\nu,\Tau}$, which corresponds to saying that the forgetful map
needs to map smooth points to smooth points (this is to be the case if
we want the map to be a map of orbifolds), can be checked by finding
an explicit description of these two open sets. Indeed, one can easily
find a list of the fixed points of the action of $\Sigma$
(resp. $\Tau$) on $Q^{\mu}$ (resp. $Q^{\nu}$), using the technique of
Lemma~8.3.2 of~\cite{delignemostowbook}.

We do not go through the trouble of writing down a general combinatorial
version of this condition, because the divisibility conditions stated
in the next section turn out to be enough to prove the result of
Theorem~\ref{thm:main}. Indeed, at least for higher-dimensional
targets, there are very few cases where the divisibility conditions
hold (see Proposition~\ref{prop:divisibility}). For each case where
the divisibility conditions do hold, we shall check whether we have a
well-defined map $Q'^{\mu,\Sigma}\rightarrow Q'^{\nu,\Tau}$.

\begin{rk}\label{rk:trivial}
  The particular case of forgetful maps between INT examples
  corresponds to the case when $\Sigma$ and $\Tau$ are both
  trivial. In that case, the condition stated in
  Lemma~\ref{lem:compat} is of course always trivially satisfied, and
  $Q'^{\mu,\Sigma}=Q^\mu$ and $Q'^{\nu,\Tau}=Q^\nu$.
\end{rk}

\subsection{Lifting to monodromy covers} \label{sec:lifttocovers}

Recall from section~\ref{sec:halfint} that, in certain cases, the map
$\B^m\rightarrow Q^\mu_{st}/\Sigma$ ramifies over points of
$Q^\mu/\Sigma$. In order to get an unramified covering, one needs to
get rid of the fixed points in $Q^{\mu}$ of the action of $\Sigma$ and
work with the open set $Q'^{\mu,\Sigma}/\Sigma\subset Q^\mu/\Sigma$
instead (see Definition~\ref{dfn:freeaction}).  In most cases,
$\pi_1(Q'^{\mu,\Sigma}/\Sigma)\simeq \pi_1(Q^{\mu}/\Sigma)$, but it
can happen that some codimension one component of the fixed point of
some $\sigma\in \Sigma$ is contained in $Q^{\mu}$, see
Remark~\ref{rk:fixedpoints}(1).

In any case, in order to lift the map to monodromy covers, we require
that $Q'^{\mu,\Sigma}$ be mapped into $Q'^{\nu,\Tau}$, and denote by
$f_*:\pi_1(Q'^{\mu,\Sigma}) \rightarrow\pi_1(Q'^{\nu,\Tau})$ the
induced map on the level of fundamental
groups. From the discussion in section~\ref{sec:halfint}, we have:
\begin{prop}\label{prop:lifttocovers}
  Suppose that $f(Q'^{\mu,\Sigma})\subset Q'^{\nu,\Tau}$, and that the
  symmetries are compatible in the sense of Lemma~\ref{lem:compat}.
  Then the map lifts to a map
  $\tilde{f}:\widetilde{Q'^{\mu,\Sigma}/\Sigma}\subset
  \widetilde{Q^\mu/\Sigma}$ if and only $f_*(K_{\mu,\Sigma})\subset
  K_{\nu,\Tau}$.
\end{prop}

Moreover, Lemma~\ref{lem:kernelhalfint} gives an explicit way to check
the condition $$f_*(K_{\mu,\Sigma})\subset K_{\nu,\Tau},$$ since it
reduces to verifying a divisibility condition between the orbifold
weights of the source and the target.

Specifically, in order to get a map of orbifolds $X_\mu\rightarrow
X_\nu$, we need to require that the weights in the target divide the
weights in the source, whenever a codimension one fixed point set of
elements of $\Sigma$ in $Q^\mu$ gets mapped onto a codimension one
fixed point set of elements of $\Tau$ in $Q^\nu$.

In other words, whenever $i\neq j$, $i,j\leq n+3$, and
$\mu_i+\mu_j<1$, we require that
\begin{equation}
  \label{eq:divisibility}
  d_{ij}^{\nu,\Tau} \textrm{ divides } d_{ij}^{\mu,\Sigma}.
\end{equation}
Recall from equation~\eqref{eq:sigmaintweights} that
condition~\eqref{eq:divisibility} means that:
\begin{itemize}
\item $(1-\nu_i-\nu_j)^{-1}$ divides $(1-\mu_i-\mu_j)^{-1}$ if $i$
      and $j$ are not in the same $\Sigma$ orbit nor in the same
      $\Tau$-orbit;
\item $(1-\nu_i-\nu_j)^{-1}$ divides $2(1-\mu_i-\mu_j)^{-1}$ if
      $i$ and $j$ are in the same $\Tau$-orbit but not in the same
      $\Sigma$-orbit;
    \item $2(1-\nu_i-\nu_j)^{-1}$ divides $2(1-\mu_i-\mu_j)^{-1}$ if
      $i$ and $j$ are in the same $\Sigma$ orbit and in the same
      $\Tau$-orbit.
\end{itemize}

\subsection{Combinatorial check}

The proof of part~(1) of Theorem~\ref{thm:main} amounts to a
combinatorial check on the list of tuples of weights that satisfy
condition $\frac{1}{2}$INT, using the results of the next few sections.

Specifically, the necessary conditions of Lemma~\ref{lem:compat} and
Proposition~\ref{prop:lifttocovers} can be checked using a computer,
so one can easily find the list of examples of forgetful maps between
any finite list of Deligne-Mostow moduli spaces.

Recall that the list of $n$-tuples satisfying condition
$\frac{1}{2}$INT is finite for $n\geq 5$ (it can be found
in~\cite{thurstonshapes}, for instance), but there are infinitely many
$4$-tuples that satisfy $\frac{1}{2}$INT (corresponding to the fact
that there are inifinitely many hyperbolic triangle groups).

However, from the translation in terms of divisibility conditions of
Proposition~\ref{prop:lifttocovers} (see the end of
section~\ref{sec:lifttocovers}), one easily gets a bound on the least
common denominator of the target weights, so in order to get maps to
strictly smaller dimension, we need only consider finitely many
$1$-dimensional moduli spaces as targets. Specifically, since all
$\frac{1}{2}$INT $k$-tuples with $k\geq 5$ have least common
denominator $\leq 42$, in order to get maps from dimension $m$ to
dimension $n$ with $m>n=1$, we need only consider $4$-tuples with
denominators $\leq 84$.

\begin{prop}\label{prop:divisibility}
  Suppose $\mu,\Sigma$ and $\nu,\Tau$ satisfy the compatibility
  conditions of Lemma~\ref{lem:compat} and the divisibility
  conditions~\eqref{eq:divisibility}. If $X_\mu$ is cocompact and
  $m\geq 3$, then up to obvious commensurabilities and symmetry, the
  pair $\mu,\nu$ is one of the following:
  $$
  \mu=(3,3,3,3,3,1)/8\textrm{ and }\nu=(3,3,3,7)/8,\quad \Sigma=\Tau=\{Id\};
  $$
  $$
  \mu=(3,3,3,3,3,1)/8\textrm{ and }\nu=(5,5,5,1)/8,\quad \Sigma=\Tau=\{Id\};
  $$
  $$
  \mu=(3,3,3,3,6,2)/10\textrm{ and }\nu=(3,3,3,3,8)/10,\quad \Sigma=\Tau=S_4.
  $$
\end{prop}

Note that the first two cases are essentially identical, since the
groups for $(3,3,3,7)/8$ and $(5,5,5,1)/8$ are isomorphic because the
sum of their respective weights is one (see Remark~\ref{rk:duality}).

Proposition~\ref{prop:divisibility} is proved by direct case by case
verification (reduced to a finite problem because of the discussion of
the beginning of this section). The author did this by writing a
computer program that generates all Deligne-Mostow sets of weights, as
well as their permutations (or rather all essential permutations,
meaning that we take the symmetry of the weights into account). For a
given pair $\mu$, $\nu$, it is of course straightforward to check the
compatibility and divisibility conditions
Proposition~\ref{prop:divisibility}.  Computer code that produces this
list is available on the author's webpage, see~\cite{derauxtautoweb}.

We now finish the proof of Theorem~\ref{thm:main}. First note that the
forgetful map corresponding to $\mu=(3,3,3,3,6,2)/10$,
$\nu=(3,3,3,3,8)/10$, with $\Sigma=\Tau=S_4$ does not yield a map of
orbifolds, because $Q'^{\mu,\Sigma}$ does not map to
$Q'^{\nu,\Tau}$. More specifically, the fixed point sets in $Q^\mu$
and $Q^\nu$ of bitranspositions have codimension one in $Q^\nu$ only,
not in $Q^\mu$ (see Lemma~8.3.2 of~\cite{delignemostowbook}, and also
Remark~\ref{rk:fixedpoints}(1)). If the map were a map of orbifolds,
the induced homomorphism between the orbifold fundamental groups would
have to map the trivial element to a nontrivial one.

Finally we need to show that the map $X_\mu\rightarrow X_\nu$ for
$\mu=(3,3,3,3,3,1)/8$, $\nu=(3,3,3,7)/8$ (and $\Sigma=\Tau=\{Id\}$)
does yield a map of orbifolds.  In that case we clearly have a map
$Q^\mu/\Sigma\rightarrow Q^\nu/\Tau$, that lifts to monodromy covers
because the divisibility condition is satisfied (see
Proposition~\ref{prop:lifttocovers}). As mentioned
in~\cite{toledomaps}, this map has a holomorphic extension to the Fox
completions $Q^\mu_{st}\rightarrow Q^\nu_{st}$ because of Hartog's
theorem, and this extension is still equivariant, so it induces a map
of orbifolds $X_\mu\rightarrow X_\nu$.

We shall give a more concrete description of the corresponding map
$X_\mu\rightarrow X_\nu$ in section~\ref{sec:example}.

\subsection{An example} \label{sec:example}

We now give some detail on the example that appears in
Theorem~\ref{thm:main}, part~(\ref{it:3dim}).
Consider
$$
\mu=(3,3,3,3,3,1)/8
$$
and 
$$
\nu=(3,3,3,7)/8.
$$
We shall choose the symmetry groups $\Sigma$, $\Tau$ to be trivial
(there are several maps obtained for various non-trivial symmetry
groups, but as mentioned above, the corresponding groups are
commensurable, see the discussion in the end of
section~\ref{sec:halfint}). Accordingly, we use the notation of
section~\ref{sec:int} rather than of section~\ref{sec:halfint}.

Note that $Q^{\nu}_{st}$ is simply a copy of $P^1$, with three
orbifold points (see~\cite{delignemostow}, p.~29). The isomorphism is
provided simply by the cross ratio
\begin{equation}\label{eq:crossratio}
  \frac{x_3-x_2}{x_3-x_1}\cdot\frac{x_4-x_1}{x_4-x_2},
\end{equation}
which is well-defined as long as no triple of points in
$\{x_1,x_2,x_3,x_4\}$ coincide (this never happens in $M^\nu_{st}$,
since any three weights of $\nu$ add up to more that one).

Now the map that sends $(x_1,\dots,x_6)$ to the cross-ratio of
$x_1,x_2,x_3,x_4$, see equation~(\ref{eq:crossratio}), is well-defined
on $M^{\mu}_{st}$ (again because no triple of weights of $\mu$ add up
to more that one), and clearly descends to $Q^\mu_{st}$.

The following result is contained in Theorem~4.1 of~\cite{klw}.
\begin{lem}
  $Q_{st}^\mu$ is a $\P^1$-bundle over $\widehat{\P^2}$, where
  $\widehat{\P^2}$ denotes $\P^2$ blown-up at a generic quadruple of
  points.
\end{lem}

\begin{pf}
  Consider the map forgetting $x_6$ from $M^{\mu}_{st}$ into
  $(\P^1)^5$. Since the stability condition allows any pair among the
  first five points to coalesce, $M^{\mu}_{st}$ maps into
  $M^{\lambda}_{st}$, where $\lambda=(2,2,2,2,2)/5$. The corresponding
  quotient $Q^{\lambda}_{st}$ is known to be $\widehat{\P^2}$ (which
  is isomorphic to $\P^1\times\P^1$ blown up in three distinct points
  of the diagonal), see~\cite{delignemostow}, Example~1, p.~33 for
  instance. 

  It is easy to check that the fibers of the corresponding map
  $M^{\mu}_{st}\rightarrow M^{\lambda}_{st}$ are all projective lines,
  and the stabilizer of a fiber (for the diagonal action of $PGL_2$)
  is trivial, so the fibers of the induced map
  $Q^{\mu}_{st}\rightarrow Q^{\lambda}_{st}$ are also projective
  lines.
\end{pf}

\begin{rk}
  $\widehat{\P^2}$ is in fact homeomorphic to $Q_{st}^\theta$ for
  various $5$-tuples $\theta$, but the corresponding map $Q_{st}^\mu
  \rightarrow Q_{st}^\theta$ is never a map of orbifolds (see
  Theorem~\ref{thm:main}).
\end{rk}

Now there is an obvious map $\widehat{\P^2}\rightarrow \P^1$, coming
from the fibration of $\P^2\setminus\{p\}$ over the $\P^1$ of lines
through $p$ (the former map contracts three of the exceptional
divisors of $\widehat{\P^2}$, and maps onto the other exceptional
divisor).

We claim that the composition $Q_{st}^\mu\rightarrow
\widehat{\P^2}\rightarrow \P^1$ can be made into a map of orbifolds,
where $\P^1$ is the orbifold ball quotient $Q_{st}^\nu$. We denote by
$f:Q_{st}^\mu\rightarrow Q_{st}^\nu$ the corresponding holomorphic
map.
\begin{figure}
  \centering \epsfig{figure=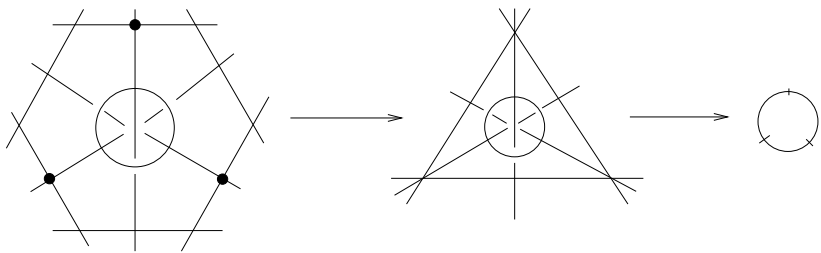,width=0.7\textwidth}
  \caption{Na{\"i}ve description of the map $\widehat{\P^2}\rightarrow\P^1$.
}
\label{fig:blowdown}
\end{figure}
The set of points of $Q_{st}^\mu$ where $df$ is not surjective
consists of three $\P^1$'s, which are the fibers over three points of
$\widehat{\P^2}$, represented by the solid dots on
Figure~\ref{fig:blowdown} (left).

Note also that the preimage of each of the three singular points on
$Q_{st}^\nu$ consists of two divisors, for instance $D_{12}$ and
$D_{34}$ map to the same point in $Q_{st}^\nu$.

We close this section by noting that the divisibility
conditions~\eqref{eq:divisibility} are trivially satisfied in this
example (recall that we take $\Sigma$ and $\Tau$ to be trivial, so
$K_\mu$ is the same as $K_{\mu,\Sigma}$, etc).  Indeed, recall that
$K_\mu$ is the normal subgroup generated by the
$\gamma_{ij}^{d_{ij}}$, where $d_{ij}$ is either $2$ or $4$. Note that
$d_{ij}=2$ only when one of $i$ or $j$ is equal to $6$, but then the
loop $\gamma_{ij}$ has trivial image. Since all the
$(1-\nu_i-\nu_j)^{-1}$ are $\pm 4$, $K_\mu$ maps into $K_\nu$.

\begin{rk}\label{rk:singular}
  Note that the divisors $D_{i5}$, $i=1,\dots,4$, $D_{i6}$,
  $i=1,\dots,5$, surject onto $Q^\nu_{st}\simeq P^1$, whereas the other
  divisors get mapped to the three orbifold points. For instance, the
  fiber over $D^\nu_{12}\in\P^1$, which corresponds to $x_1=x_2$ in
  $Q^\nu_{st}$, is given by the union of the divisors $D^\mu_{12}$ and
  $D^\mu_{34}$ in $Q^\mu_{st}$ (these are both projective planes and their
  intersection is a projective line).
\end{rk}

\section{Maps between non-compact examples}

The same construction works for non compact moduli spaces, and one
gets maps that have the same behavior as the ones between compact
ones. Here the divisibility condition is easily adapted to allow for
infinite weights if some pairs of weights add up to exactly one.  The
list of pairs of weights with symmetry (at least one of which is
non-compact, with $m>2$ or $n>1$) that satisfy the compatibility and
divisibility conditions is the following:
\begin{itemize}
\item $\mu=(2,2,2,3,3)/6$, $\Sigma=S_3$\\
  $\nu=(2,2,2,1,5)/6$, $\Tau=S_3$;
\item $\mu=(4,4,4,5,7)/12$, $\Sigma=S_3$\\
  $\nu=(2,2,2,1,5)/6$, $\Tau=S_3$;
\item $\mu=(2,2,3,3,1,1)/6$, $\Sigma=S_2$\\
  $\nu=(1,7,7,9)/12$, $\Tau=\{id\}$;
\item $\mu=(2,2,3,3,1,1)/6$, $\Sigma=S_2$\\
  $\nu=(1,3,4,4)/6$, $\Tau=\{id\}$.
\end{itemize}

Note that in particular one does not get more pairs $(m,n)$ of
dimensions that are related by a surjective map of orbifolds coming
from a forgetful map than in the compact case; in other words part~(i)
of the statement of Theorem~\ref{thm:main} remains true for non-compact
examples.

\bibliographystyle{amsplain}

\begin{thebibliography}{10}

\bibitem{bergeron}
N.~Bergeron, F.~Haglund, and D.~Wise, \emph{Hyperplane sections in arithmetic
  hyperbolic manifolds}, Preprint, 2008.

\bibitem{delignemostow}
P.~Deligne and G.~D. Mostow, \emph{Monodromy of hypergeometric functions and
  non-lattice integral monodromy}, Inst. Hautes {\'E}tudes Sci. Publ. Math.
  \textbf{63} (1986), 5--89.

\bibitem{delignemostowbook}
\bysame, \emph{Commensurabilities among lattices in {P}{U}(1,n)}, Annals of
  Mathematics Studies, vol. 132, Princeton Univ. Press, Princeton, 1993.

\bibitem{derauxtautoweb}
M.~Deraux, \emph{web page, http://www.math.utah.edu/{$\sim$}deraux/java}.

\bibitem{derauxthesis}
\bysame, \emph{Complex surfaces of negative curvature}, Ph.D. thesis, Univ. of
  Utah, 2001.

\bibitem{derauxmathann}
\bysame, \emph{On the universal cover of certain exotic {K}{\"a}hler surfaces
  of negative curvature}, Math. Ann. \textbf{329} (2004), no.~4, 653--683.

\bibitem{deraux3fold}
\bysame, \emph{A negatively curved {K}\"ahler threefold not covered by the
  ball}, Inv. Math. \textbf{160} (2005), no.~3, 501--525.

\bibitem{klw}
F.~C. Kirwan, R.~Lee, and S.~H. Weintraub, \emph{Quotients of the complex ball
  by discrete groups}, Pacific J. Math. \textbf{130} (1987), 115--141.

\bibitem{koziarzmok}
V.~Koziarz and N.~Mok, \emph{Nonexistence of holomorphic submersions between
  complex unit balls equivariant with respect to a lattice and their
  generalizations}, Preprint, \verb|math.AG/08042122|, 2008.

\bibitem{livne}
R.~A. Livn{\'e}, \emph{On certain covers of the universal elliptic curve},
  Ph.D. thesis, Harvard University, 1981.

\bibitem{looijenga}
E.~Looijenga, \emph{Uniformization by {L}auricella functions---an overview of
  the theory of {D}eligne-{M}ostow}, Arithmetic and geometry around
  hypergeometric functions (R.-P. Holzapfel, A.~M. Ulud{\u a}g, and M.~Yoshida,
  eds.), Progress in Mathematics, vol. 260, Birkh{\"a}user, Basel, 2007,
  pp.~207--244.

\bibitem{mokannals}
N.~Mok, \emph{Uniqueness theorems of hermitian metrics of seminegative
  curvature on quotients of bounded symmetric domains}, Ann. of Math.
  \textbf{125} (1987), no.~2, 105--152.

\bibitem{mostowihes}
G.~D. Mostow, \emph{Generalized {P}icard lattices arising from half-integral
  conditions}, Inst. Hautes {\'E}tudes Sci. Publ. Math. \textbf{63} (1986),
  91--106.

\bibitem{mostowdiscontinuous}
\bysame, \emph{On discontinuous action of monodromy groups on the complex
  $n$-ball}, J. Amer. Math. Soc. \textbf{1} (1988), 555--586.

\bibitem{mostowsiu}
G.~D. Mostow and Y.~T. Siu, \emph{A compact {K\"ahler} surface of negative
  curvature not covered by the ball}, Ann. of Math. \textbf{112} (1980),
  321--360.

\bibitem{siuconjecture}
Y.-T. Siu, \emph{Some recent results in complex manifold theory related to
  vanishing theorems for the semipositive case}, Arbeitstagung Bonn 1984
  (F.~Hirzebruch, J.~Schwermer, and S.~Suter, eds.), Lecture Notes in
  Mathematics, vol. 1111, Springer Verlag, 1985, pp.~169--192.

\bibitem{thurstonshapes}
W.~P. Thurston, \emph{Shapes of polyhedra and triangulations of the sphere},
  Geometry and Topology Monographs \textbf{1} (1998), 511--549.

\bibitem{toledomaps}
D.~Toledo, \emph{Maps between complex hyperbolic surfaces}, Geom. Ded.
  \textbf{97} (2003), 115--128.

\end{thebibliography}

\providecommand{\bysame}{\leavevmode\hbox to3em{\hrulefill}\thinspace}
\providecommand{\MR}{\relax\ifhmode\unskip\space\fi MR }
\providecommand{\MRhref}[2]{%
  \href{http://www.ams.org/mathscinet-getitem?mr=#1}{#2}
}
\providecommand{\href}[2]{#2}

\end{document}